\documentclass[10pt]{amsart}
\usepackage{indentfirst,latexsym,bm}
\usepackage{amsfonts}

\usepackage{amssymb}
\usepackage{amsmath}
\usepackage{hyperref}
\usepackage{dsfont}
\usepackage{leftidx}
\usepackage{amsthm}
\usepackage{geometry}
\usepackage{amsmath,amscd}
\usepackage{tikz}
\usepackage[all,cmtip]{xy}

\usepackage{latexsym,amsmath,amssymb,amsfonts}  
\usepackage{color}
\usepackage{color,xcolor}
\usepackage{colordvi}
\usepackage{xspace}
      
\usepackage{graphicx}
\usepackage{pstricks}
\usepackage{pgf}
\usepackage{tikz}
\usepackage{xkeyval}
\usepackage{tikz}
\usetikzlibrary{mindmap,trees,shadows}

\begin{document}
 \thispagestyle{empty}

\title[Noetherian connected braided Hopf algebras are affine]
{Noetherian connected braided Hopf algebras are affine}

\author{
{Huan Jia $^{1}$}
and 
{Yinhuo Zhang $^{2}$}}

\address{$^{1}$ Department of Mathematics, Suqian University, Suqian City 223800, Jiangsu Province, China}

\address{$^{2}$ Department of Mathematics and Statistics,
University of Hasselt, Universitaire Campus, 3590 Diepenbeek, Belgium}

\address{Huan Jia $^{1}$}
\email{huan.jia@squ.edu.cn}

\address{Yinhuo Zhang $^{2}$}
\email{yinhuo.zhang@uhasselt.be}
 
\date{\today}
\maketitle

\newtheorem{theorem}{Theorem}[section]
\newtheorem{proposition}[theorem]{Proposition}
\newtheorem{lemma}[theorem]{Lemma}
\newtheorem{corollary}[theorem]{Corollary}
\theoremstyle{definition}
\newtheorem{definition}[theorem]{Definition}
\newtheorem{example}[theorem]{Example}
\newtheorem{remark}[theorem]{Remark}

\newcommand{\K}{\mathds{k}}
\newcommand{\X}{\langle X \rangle}

\newcommand{\id}{\operatorname{id}}
\newcommand{\lex}{\operatorname{lex}}
\newcommand{\rank}{\operatorname{rank}}
\newcommand{\lOr}{\langle O \rangle}
\newcommand{\rf}{\operatorname{rf}}
\newcommand{\HYD}{{}_{H}^{H}\mathcal{YD}}

\theoremstyle{plain}
\newcounter{maint}
\renewcommand{\themaint}{\Alph{maint}}
\newtheorem{mainthm}[maint]{Theorem}

\theoremstyle{plain}
\newtheorem*{proofthma}{Proof of Theorem A}
\newtheorem*{proofthmb}{Proof of Theorem B}

\begin{abstract}
The well-known Wu--Zhang question asks whether every Noetherian Hopf algebra over a field is affine (finitely generated as an algebra). In this paper, we extend this fundamental problem to the braided setting within the monoidal category ${}_{H}^{H}\mathcal{YD}$ of left Yetter--Drinfel'd modules over a Hopf algebra $H$. By developing term-rewriting methods based on reduction orderings and prime-word factorizations in free prebraided algebras, we establish explicit expansions for YD-triangular comultiplications. We prove that if $(\mathbb{K}\langle X \rangle, \epsilon)$ is an augmented free algebra in ${}_{H}^{H}\mathcal{YD}$ and $I \subseteq \ker \epsilon$ is an ideal invariant under a YD-triangular comultiplication, then every $I$-reducible letter---and consequently every non-empty word in $\langle X \rangle$---can be expressed modulo $I$ as a linear combination of products of $I$-irreducible letters strictly smaller with respect to the reduction order. 

As a primary consequence, we prove that any right Noetherian connected prebraided Hopf algebra $B$ in ${}_{H}^{H}\mathcal{YD}$ is affine, generated by the image of a finite set $X_{I}$ of $I$-irreducible letters. As an application, we show that the braided component $B = A^{\mathrm{co}\,\pi}$ of a right Noetherian Radford biproduct $A = B \sharp H$ is affine whenever $A_{(0)} \subseteq \jmath(H)$, reducing the affineness of $A$ entirely to that of $H$.

\end{abstract}

\section*{Introduction}

Let $\mathbb{K}$ be an arbitrary field. A $\mathbb{K}$-algebra is called \textit{affine} if it is finitely generated as an algebra. In 2003, Q.-S.~Wu and J.~J.~Zhang \cite{WZ2003} established that every affine Noetherian PI Hopf algebra is Artin--Schelter Gorenstein, and raised the following well-known open question for general Noetherian Hopf algebras (see also \cite{A2023, Br2007, Br2026, BG2014, G2013}):

\begin{quote}
\textbf{Question} (\cite[Question~0.1]{WZ2003}).
\textit{Is every Noetherian Hopf algebra over $\mathbb{K}$ affine?}
\end{quote}

Molnar \cite{M1975} previously established that every commutative Noetherian Hopf algebra over a field is affine (with the cocommutative case following dually). In 2017, Goodearl and Zhang \cite{GZ2017} proved that a locally affine Noetherian Hopf algebra is affine provided it is faithfully flat over a Hopf subalgebra, establishing in particular that every Noetherian pointed Hopf algebra is affine.

Recently, Jia and Zhang \cite{JZ2025-1} proved that a Noetherian algebra-graded Hopf algebra is affine if and only if its degree-zero component is affine, and deduced that the braided component of a Noetherian graded Hopf algebra is affine. To address the general (non-graded) setting of the Wu--Zhang question, Jia and Zhang \cite{JZ2025-2, JZ2026} introduced term-rewriting tools based on (mirror) reduction orderings and prime word factorizations. Combined with the generalized lifting method \cite{AC2013}, they proved in \cite{JZ2026} that a Noetherian Hopf algebra is affine if and only if its Hopf coradical is affine, thereby reducing the Wu--Zhang question to the first step of the standard filtration.

Building on these developments in the ordinary setting, it is natural to extend the Wu--Zhang question to the braided framework: \textit{Is every Noetherian braided Hopf algebra affine?} Beyond its intrinsic algebraic interest, this question plays a critical structural role because many ordinary Hopf algebras arise as deformations or liftings of Radford biproducts (bosonizations) constructed from connected braided Hopf algebras—most notably Nichols algebras, pre-Nichols algebras, and braided tensor algebras (see, e.g., \cite{AC2013, AS1998}).

In this paper, we develop reduction-ordering and prime-word rewriting techniques within the prebraided monoidal category ${}_{H}^{H}\mathcal{YD}$ to address this problem. Our main result establishes that any right Noetherian connected prebraided (and in particular, braided) Hopf algebra is affine.

\subsection*{Organization of the Paper}

Let $(X, \prec_{X})$ be a well-ordered set, let $\langle X \rangle$ denote the free monoid on $X$, and let $\mathbb{K}\langle X \rangle$ denote the free algebra on $X$.

In Section~\ref{sec:1}, we review necessary preliminaries on Hopf algebras, the category ${}_{H}^{H}\mathcal{YD}$ of left Yetter--Drinfel'd $H$-modules, and Radford biproducts.

In Section~\ref{sec:2}, we adapt the notions of reduction orders and prime words from \cite{JZ2025-2, JZ2026} to our prebraided setting. We study $I$-irreducible words in $\langle X \rangle$ with respect to $\prec_{X}$ and prove that if an augmented algebra is right Noetherian, then the underlying set of $I$-irreducible letters is finite (Lemma~\ref{lem:Noetherian-irrls-finite}).

In Section~\ref{sec:3}, we apply reduction orders and prime words to determine the explicit comultiplication structure on the free algebra $\mathbb{K}\langle X \rangle$ in ${}_{H}^{H}\mathcal{YD}$ (Lemmas~\ref{lem:coproduct-omega-1} and \ref{lem:coproduct-omega-2}).

In Section~\ref{sec:4}, we utilize this comultiplication structure to establish a crucial reduction lemma for prime words (Lemma~\ref{lem:wr-sum-smaller-words}). We then demonstrate that every reducible letter can be expressed as a linear combination of products of letters of strictly smaller order (Lemma~\ref{lem:rl-sum-smaller-irrls}). As a consequence, we obtain the following primary representation theorem:

\medskip
\noindent
\textbf{Theorem A} (Theorem~\ref{thm:rl-sum-smaller-irrls})\textbf{.}
\textit{Let $(\mathbb{K}\langle X \rangle, \epsilon)$ be the augmented free algebra in ${}_{H}^{H}\mathcal{YD}$, and let $I \subseteq \ker\epsilon$ be a proper ideal of $\mathbb{K}\langle X \rangle$. Suppose $\Delta$ is a YD-triangular comultiplication on $\mathbb{K}\langle X \rangle$ such that $\Delta(I) \subseteq I \otimes \mathbb{K}\langle X \rangle + \mathbb{K}\langle X \rangle \otimes I$. If $x \in X$ is an $I$-reducible letter, then there exist $I$-irreducible letters $x_{1}, \ldots, x_{n} \in X_{I}$ and scalars $k_{x_{1}\cdots x_{n}} \in \mathbb{K}$ such that
\[
x \in \sum k_{x_{1}\cdots x_{n}} x_{1} \cdots x_{n} + I, \qquad \text{with } x_{1}, \ldots, x_{n} \prec_{X} x.
\]
Consequently, every non-empty word in $\langle X \rangle$ can be expressed modulo $I$ as a linear combination of products of $I$-irreducible letters.}
\medskip

In Section~\ref{sec:5}, we apply these foundational results to establish our main generating theorem:

\medskip
\noindent
\textbf{Theorem B} (Theorem~\ref{thm:Noetherian-connected-braided-HA-affine})\textbf{.}
\textit{Let $H$ be a Hopf algebra over $\mathbb{K}$, and let $B$ be a connected Hopf algebra in ${}_{H}^{H}\mathcal{YD}$. Then the following assertions hold:
\begin{itemize}
  \item[(a)] There exists a generating set $X$ constructed as in Theorem~\ref{thm:rl-sum-smaller-irrls} such that the set $\pi(X_{I})$ generates $B$ as an algebra, where $\pi \colon \mathbb{K}\langle X \rangle \to B$ is the canonical algebra projection and $I := \ker\pi$.
  \item[(b)] If $B$ is right Noetherian, then $B$ is affine.
\end{itemize}}
\medskip

As an immediate structural application, we deduce that the braided component of a Noetherian Radford biproduct is affine under a natural coradical inclusion hypothesis (Corollary~\ref{cor:Noetherian-Radford-biproduct-affine}).

\medskip
\noindent
Several conceptual techniques utilized in this paper are inspired by foundational frameworks in non-commutative term-rewriting and PBW theories \cite{Ji2023, K1999, Ro1999, U2004, Z2024, ZSL2020}, particularly the notion of irreducible words under graded lexicographic orders and the use of free (braided) comultiplications in establishing PBW bases.

\section{Lexicographic Order and Prebraided Hopf Algebras}\label{sec:1}

\subsection{Lexicographic Order on Words}

For a set $X$, let $\langle X \rangle$, $\langle X \rangle^{+}$, and $\mathbb{K}\langle X \rangle$ denote the free monoid, the free semigroup, and the free algebra generated by $X$, respectively. The empty word in $\langle X \rangle$ is denoted by $1$. Clearly, $\langle X \rangle = \langle X \rangle^{+} \cup \{1\}$, and $\mathbb{K}\langle X \rangle^{+}$ is the subspace of $\mathbb{K}\langle X \rangle$ spanned by $\langle X \rangle^{+}$.

A word $\beta \in \langle X \rangle$ is called a \textit{factor} of $\alpha \in \langle X \rangle$ if there exist $\gamma, \eta \in \langle X \rangle$ such that $\alpha = \gamma\beta\eta$. The factor $\beta$ is a \textit{prefix} (resp.\ \textit{suffix}) of $\alpha$ if $\gamma = 1$ (resp.\ $\eta = 1$). If $\beta \neq \alpha$, then $\beta$ is called a \textit{proper} factor (resp.\ proper prefix, proper suffix) of $\alpha$.

For a word $\alpha = x_{1} \cdots x_{n} \in \langle X \rangle$ with $x_{1}, \ldots, x_{n} \in X$, we define:
\begin{itemize}
  \item $\alpha(y)$ as the number of occurrences of a letter $y \in X$ in $\alpha$, i.e., $\alpha(y) := \# \{i \mid x_{i} = y,\ 1 \le i \le n\}$;
  \item $l(\alpha)$ as the leftmost letter of $\alpha$, i.e., $l(\alpha) := x_{1}$;
  \item $|\alpha|$ as the \textit{length} of $\alpha$, i.e., $|\alpha| := n$.
\end{itemize}
For each integer $m \ge 1$, we write $X^{m} := \{\alpha \in \langle X \rangle \mid |\alpha| = m\}$.

Let $(X, \prec_{X})$ be a well-ordered set. The \textit{lexicographic order} $\prec_{\mathrm{lex}_{X}}$ on $\langle X \rangle$ induced by $\prec_{X}$ is defined as follows (see \cite{L1997}). For any two words $\alpha, \beta \in \langle X \rangle$,
\begin{align*}
\alpha \prec_{\mathrm{lex}_{X}} \beta \Longleftrightarrow 
\begin{cases}
\alpha \text{ is a proper prefix of } \beta, \text{ or} \\
\alpha = \gamma x \eta \text{ and } \beta = \gamma y \zeta \text{ with } x, y \in X, x \prec_{X} y, \text{ and } \gamma, \eta, \zeta \in \langle X \rangle.
\end{cases}
\end{align*}

\subsection{Prebraided Bialgebras and Hopf Algebras}

We refer to \cite{Ra2012, S1969} for standard background on Hopf algebras and prebraided Hopf algebras. The category of $\mathbb{K}$-vector spaces is denoted by $\mathbb{K}\text{-Vec}$.

Let $(C, \Delta_{C}, \epsilon_{C})$ be a coalgebra over $\mathbb{K}$. We use Sweedler's notation for the coproduct, writing $\Delta_{C}(c) = \sum c_{(1)} \otimes c_{(2)}$ for $c \in C$. The \textit{coradical} of $C$, denoted by $C_{(0)}$, is the sum of all simple subcoalgebras of $C$. The coalgebra $C$ is called \textit{connected} if $C_{(0)} = \mathbb{K}$. The \textit{coradical filtration} $\{C_{(n)}\}_{n=0}^{\infty}$ of $C$ is defined inductively for $n \ge 1$ by 
\[
C_{(n)} := C_{(n-1)} \wedge C_{(0)} = \Delta_{C}^{-1}(C_{(n-1)} \otimes C + C \otimes C_{(0)}).
\]
This forms an ascending coalgebra filtration of $C$.

Let $H$ be a Hopf algebra over $\mathbb{K}$ with antipode $S_{H}$. We denote by ${}_{H}^{H}\mathcal{YD}$ the category of left Yetter--Drinfel'd $H$-modules. Its objects are triples $(M, \cdot, \rho)$, where $(M, \cdot)$ is a left $H$-module and $(M, \rho)$ is a left $H$-comodule with coaction $\rho(m) = \sum m_{(-1)} \otimes m_{(0)}$, satisfying the compatibility condition
\begin{align*}
\sum h_{(1)} m_{(-1)} \otimes h_{(2)} \cdot m_{(0)} = \sum (h_{(1)} \cdot m)_{(-1)} h_{(2)} \otimes (h_{(1)} \cdot m)_{(0)}
\end{align*}
for all $h \in H$ and $m \in M$.

For any objects $M, N \in {}_{H}^{H}\mathcal{YD}$, the tensor product $M \otimes N$ becomes an object in ${}_{H}^{H}\mathcal{YD}$ via the module action and coaction
\begin{align*}
h \cdot (m \otimes n) &= \sum h_{(1)} \cdot m \otimes h_{(2)} \cdot n, \\
\rho_{M \otimes N}(m \otimes n) &= \sum m_{(-1)} n_{(-1)} \otimes m_{(0)} \otimes n_{(0)},
\end{align*}
for $h \in H$, $m \in M$, and $n \in N$. The category ${}_{H}^{H}\mathcal{YD}$ is prebraided monoidal, with prebraiding $c_{M,N} \colon M \otimes N \to N \otimes M$ given by
\begin{align*}
c_{M,N}(m \otimes n) = \sum m_{(-1)} \cdot n \otimes m_{(0)}, \quad m \in M, ~n \in N.
\end{align*}
In particular, if the antipode $S_{H}$ is bijective, then $c_{M,N}$ is an isomorphism for all $M, N$, and ${}_{H}^{H}\mathcal{YD}$ is a braided monoidal category.

An algebra in ${_{H}^{H}\mathcal{YD}}$ is a triple $(A, m_{A}, \eta_{A})$, where $A$ is a left $H$-module algebra and a left $H$-comodule algebra. 
Given algebras $(A,m_{A},\eta_{A})$ and $(B,m_{B},\eta_{B})$ in ${_{H}^{H}\mathcal{YD}}$. 
The tensor product $(A \otimes B, m_{A \otimes B}, \eta_{A\otimes B})$ is an algebra in ${_{H}^{H}\mathcal{YD}}$, whose algebra structure is defined by 
\begin{align*}
m_{A \otimes B} &= (m_{A} \otimes m_{B})\circ (\operatorname{id} \otimes c_{A,B} \otimes \operatorname{id}), \qquad \eta_{A \otimes B} = \eta_{A} \otimes \eta_{B},
\end{align*}
that is, for $a,c\in A$ and $b,d \in B$,
\begin{align*}
(a \otimes b)(c \otimes d) &= a(b_{(-1)} \cdot c) \otimes b_{(0)}d, \qquad \eta_{A\otimes B}(1) = 1 \otimes 1.
\end{align*}
Sometimes we write $A \otimes^{c} B$ for the tensor products of algebras $A,B$ in ${_{H}^{H}\mathcal{YD}}$ to distinguish it from the tensor product in $\mathbb{K}\text{-Vec}$.

\vspace{1em}

A coalgebra in ${_{H}^{H}\mathcal{YD}}$ is a triple $(C, \Delta_{C}, \epsilon_{C})$, where $C$ is a left $H$-module coalgebra and a left $H$-comodule coalgebra. Given coalgebras $(C,\Delta_{C},\epsilon_{C})$ and $(D,\Delta_{D},\epsilon_{D})$ in ${_{H}^{H}\mathcal{YD}}$. The tensor product $(C \otimes D, \Delta_{C \otimes D}, \epsilon_{C\otimes D})$ is a coalgebra in ${_{H}^{H}\mathcal{YD}}$, whose coalgebra structure is defined by 
\begin{align*} 
\Delta_{C \otimes D} &= (\operatorname{id} \otimes c_{C,D} \otimes \operatorname{id}) \circ (\Delta_{C} \otimes \Delta_{D}), \qquad \epsilon_{C \otimes D} = \epsilon_{C} \otimes \epsilon_{D},
\end{align*}
that is, for $c\in C$ and $d\in D$,
\begin{align*}
\Delta_{C \otimes D}(c \otimes d) &= (c_{(1)} \otimes c_{(2)(-1)} \cdot d_{(1)}) \otimes (c_{(2)(0)} \otimes d_{(2)}), \qquad \epsilon_{C\otimes D}(c \otimes d) = \epsilon_{C}(c) \epsilon_{D}(d). 
\end{align*}
Sometimes we write $C \otimes^{c} D$ for the tensor product of coalgebras $C,D$ in ${_{H}^{H}\mathcal{YD}}$.

\vspace{1em}

A bialgebra in ${_{H}^{H}\mathcal{YD}}$ is a tuple $(B, m_{B}, \eta_{B}, \Delta_{B}, \epsilon_{B})$, where $(B, m_{B}, \eta_{B})$ is an algebra in ${_{H}^{H}\mathcal{YD}}$, and $(B, \Delta_{B}, \epsilon_{B})$ is a coalgebra in ${_{H}^{H}\mathcal{YD}}$, satisfying the compatibility conditions:
\begin{align*}
\Delta_{B}(1) &= 1 \otimes 1, \quad \Delta_{B}(ab) = a_{(1)}(a_{(2)(-1)} \cdot b_{(1)}) \otimes a_{(2)(0)} b_{(2)}, \\
\epsilon_{B}(1) &= 1, \quad \qquad \epsilon_{B}(ab) = \epsilon_{B}(a)\epsilon_{B}(b), \qquad a,b\in B.
\end{align*}

A Hopf algebra in ${_{H}^{H}\mathcal{YD}}$ is a tuple $(B, m_{B}, \eta_{B}, \Delta_{B}, \epsilon_{B}, S_{B})$, where $(B, m_{B}, \eta_{B}, \Delta_{B}, \epsilon_{B})$ is a bialgebra in ${_{H}^{H}\mathcal{YD}}$, and the map $S_{B}: B \rightarrow B$ (called the antipode of $B$) satisfies:
\begin{align*}
\sum S_{B}(b_{(1)})b_{(2)} = \epsilon_{B}(b) 1 = \sum b_{(1)} S_{B}(b_{(2)}), \quad b \in B.
\end{align*}

Let $B$ be a bialgebra (resp. Hopf algebra) in ${_{H}^{H}\mathcal{YD}}$. The Radford biproduct $B \sharp H := B \otimes H$ of $B$ by $H$ is a bialgebra (resp. Hopf algebra) over $\mathbb{K}$, with the smash product and the smash coproduct defined by: for $a, b \in B$ and $g, h \in H$,
\begin{align*}
(a \sharp g)(b \sharp h) = a (g_{(1)} \cdot b) \sharp g_{(2)} h, \qquad 
\Delta_{B\sharp H} (b \sharp h) = (b_{(1)} \sharp b_{(2)(-1)} h_{(1)}) \otimes (b_{(2)(0)} \sharp h_{(2)}).
\end{align*}

Conversely, many bialgebras (resp. Hopf algebras) are the Radford biproducts of a Hopf subalgebra and a bialgebra (resp. Hopf algebra) in ${_{H}^{H}\mathcal{YD}}$.

\begin{lemma}\cite{Ra1985,Ra2012}\label{lem:Radfordbiproduct}
Let $A$ be a bialgebra and $H$ a Hopf algebra over $\mathbb{K}$. Suppose that $\pi: A \rightarrow H$ and $\jmath: H \rightarrow A$ are bialgebra maps such that $\pi \circ \jmath = \operatorname{id}_{H}$. Then there exists a bialgebra $(B, \cdot, \rho) \in {_{H}^{H}\mathcal{YD}}$ such that as bialgebras,
\begin{align*}
B \sharp H &\cong A \quad (b \otimes h \mapsto b\jmath(h)).
\end{align*}
Here, $B := \Pi(A) = A^{\operatorname{co}\,\pi}$ via the coalgebra map $\Pi = \operatorname{id} * (\jmath \circ S \circ \pi)$, and its structure is given by
\begin{align*}
\Delta_{B} &= (\Pi \otimes \operatorname{id}) \circ \left.\Delta_{A}\right|_B, \\
h \cdot b &= \jmath(h_{(1)})b\jmath(S_{H}(h_{(2)})), \\
\rho(b) &= \pi(b_{(1)}) \otimes b_{(2)}, \qquad h \in H,\ b \in B.
\end{align*}
\end{lemma}

The following lemmas collect some simple observations which will be used in the sequel.

\begin{lemma}\label{lem:B_(n)-algebra-filtration}
Let $H$ be a Hopf algebra over $\mathbb{K}$, and let $B$ be a Hopf algebra in ${_{H}^{H}\mathcal{YD}}$. If $B_{(0)}$ is a left $H$-subcomodule (resp. $H$-submodule, Yetter--Drinfel'd $H$-submodule) of $B$, then $B_{(n)}$ is a left $H$-subcomodule (resp. $H$-submodule, Yetter--Drinfel'd $H$-submodule) of $B$. Moreover, $B_{(0)}$ is a sub-bialgebra of $B$ if and only if the coradical filtration $\{B_{(n)}\}_{n\ge 0}$ is a bialgebra filtration of $B$ within ${_{H}^{H}\mathcal{YD}}$.
\end{lemma}

\begin{lemma}\label{lem:B-Noetherian}
Let $H$ be a bialgebra over $\mathbb{K}$, and let $B$ be a left $H$-module algebra. If the smash product $B \sharp H$ is right Noetherian, then $B$ is right Noetherian. 
\end{lemma}

\section{Reduction Order, Prime Words, and Irreducible Letters}\label{sec:2}

For our purposes, we adapt the notions of reduction order, reduction-factorization, and prime words introduced in \cite{JZ2025-2, JZ2026}.

Let $O$ be a set. For each $z \in O$, assign a set $C_z$ of formal symbols such that $z \in C_z$. Suppose these sets $C_z$ are pairwise disjoint, and define
\begin{align*}
X := \bigcup_{z \in O} C_z.
\end{align*}
We define a map $C: X \rightarrow O$ by 
\begin{align*}
C(x) = z, \qquad z \in O,\ x \in C_{z}.
\end{align*}
We refer to $O$ as the \textit{original set} of $X$, and call $z$ the \textit{original element} of $x$.
The map $C$ naturally lifts to a monoid map $C: \X \rightarrow \lOr$.

Given a well-order $\prec_{O}$ on $O$. For each $z \in O$, we equip $C_z$ with a well-order $\prec_{C_z}$ such that $z$ is the least element. 
Define an order $\prec_{X}$ on $X$ as follows. For letters $x, y \in X$,
\begin{align*}
x \prec_{X} y 
\Longleftrightarrow 
\begin{cases}
C(x) \prec_{O} C(y), \ \  \text{or}  &\\
C(x) = C(y) = z \text{ and } x \prec_{C_z} y.
\end{cases}
\end{align*}

We now introduce the following reduction orders:

\begin{definition}\label{def:reduction-order}
The \textit{reduction order} $\prec_{\lOr}$ on $\lOr$ is defined as follows. For any non-empty words $\alpha, \beta \in \lOr$, we set $1 \prec_{\lOr} \alpha$, and
\begin{align*}
\alpha \prec_{\lOr} \beta \Longleftrightarrow 
\begin{cases}
l(\alpha) \prec_{O} l(\beta), \text{ or }  \\
l(\alpha) = l(\beta) \text{ and } 
\begin{cases}
|\alpha| < |\beta|, \text{ or }  \\
|\alpha| = |\beta| \text{ and } \alpha \prec_{\lex_{O}} \beta.
\end{cases}
\end{cases}
\end{align*}
The \textit{reduction order} $\prec_{\X}$ on $\X$ is defined as follows. For any words $\alpha, \beta \in \X$, 
\begin{align*}
\alpha \prec_{\X} \beta \Longleftrightarrow 
\begin{cases}
C(\alpha) \prec_{\lOr} C(\beta),  \text{ or } \\
C(\alpha) = C(\beta) \text{ and } \alpha \prec_{\lex_{X}} \beta.
\end{cases}
\end{align*}
\end{definition}

By definition, each word $\alpha \in \X$ satisfies $C(\alpha) \preceq_{\X} \alpha$. 

The above ordered sets are well-ordered:

\begin{lemma}\label{lem:well-orders}
$(X, \prec_{X})$, $(\lOr, \prec_{\lOr})$, and $(\X, \prec_{\X})$ are well-ordered.
\end{lemma}

The following orders are mutually compatible:

\begin{lemma}\label{lem:compatibility-orders}
The following assertions hold:
\begin{itemize}
\item[(a)] For any letters $x, y \in O$, $x \prec_{\lOr} y$ if and only if $x \prec_{O} y$.
\item[(b)] For any letters $x, y \in O$, $x \prec_{X} y$ if and only if $x \prec_{O} y$.
\item[(c)] For any letters $x, y \in X$, $x \prec_{\X} y$ if and only if $x \prec_{X} y$.
\item[(d)] For any words $\alpha, \beta \in \lOr$, $\alpha \prec_{\X} \beta$ if and only if $\alpha \prec_{\lOr} \beta$.
\end{itemize}
\end{lemma}

Applying Lemma~\ref{lem:compatibility-orders} to words yields the following immediate consequence:

\begin{corollary}\label{cor:well-orders-compatibility-2}
Let $\alpha$ be a non-empty word in $\X$ and $y$ a letter in $O$. Then the following conditions are equivalent:
\begin{itemize}
	\item[(a)] $\alpha \prec_{\X} y$;
	\item[(b)] $l(\alpha) \prec_{X} y$;
	\item[(c)] $l(C(\alpha)) = C(l(\alpha)) \prec_{O} y$;
	\item[(d)] $C(\alpha) \prec_{\lOr} y$.
\end{itemize}
\end{corollary}

By Corollary~\ref{cor:well-orders-compatibility-2}, we obtain that reduction orders are right-compatible.

\begin{lemma}\label{lem:ac<bc}
Let $\alpha, \beta$ be non-empty words in $\lOr$ (resp. $\X$) with $\alpha \prec_{\lOr} \beta$ (resp. $\alpha \prec_{\X} \beta$). Then $\alpha\gamma \prec_{\lOr} \beta\gamma$ for any $\gamma \in \lOr$ (resp. $\X$).
\end{lemma}

By virtue of Lemma~\ref{lem:compatibility-orders}, we obtain the following minimality property for letters in $O$.

\begin{lemma}\label{lem:x-least}
Let $x \in O$ and $L_{x} := \{\alpha \in \X \mid C(l(\alpha)) = x\}$. Then $x$ is the least element of the set $L_{x}$ with respect to $\prec_{\X}$. Moreover, for a letter $a \in X$, the following conditions are equivalent:
\begin{itemize}
\item[(a)] $a$ is the least non-empty word of $\X$ with respect to $\prec_{\X}$;
\item[(b)] $a$ is the least non-empty word of $\lOr$ with respect to $\prec_{\lOr}$;
\item[(c)] $a$ is the least letter of $X$ with respect to $\prec_{X}$;
\item[(d)] $a$ is the least letter of $O$ with respect to $\prec_{O}$.
\end{itemize}
\end{lemma}

We next define the following factorization of a word with respect to the reduction order $\prec_{\X}$:

\begin{definition}\label{def:reduction-factorization}
Let $\alpha$ be a non-empty word in $\X$. The \textit{reduction-factorization} of $\alpha$, denoted by $\operatorname{rf}(\alpha) = (\alpha_{L}, \alpha_{R})$, is the unique decomposition 
\[ \alpha = \alpha_{L}\alpha_{R}, \] 
where $\alpha_{R}$ is the greatest suffix of $\alpha$ with respect to the reduction order $\prec_{\X}$. 
We call the word $\alpha$ \textit{prime} if $\alpha_{R} = \alpha$.
The set of prime words in $\X$ is denoted by $\X_{p}$.
\end{definition}

By definition, each letter $x \in X$ is trivially prime because $x_{R} = x$.

For each word $\alpha \in \X$, let $m_{\alpha}$ denote the greatest letter of $\alpha$ with respect to $\prec_{X}$. 
The reduction-factorization of a word is related to the original element of its greatest letter:

\begin{lemma}\label{lem:properties-alpha_R}
Let $\alpha$ be a non-empty word in $\X$, and let $\operatorname{rf}(\alpha) = (\alpha_{L}, \alpha_{R})$. Set $m := C(m_{\alpha})$. Then the following assertions hold:
\begin{itemize}
    \item[(a)] $C(C(\alpha)) = C(\alpha)$, $C(l(\alpha)) = l(C(\alpha))$, and $C(m_{\alpha}) = m_{C(\alpha)}$.
    \item[(b)] $C(\alpha_{L})(m) = 0$ and $l(C(\alpha_{R})) = m$. 
    \item[(c)] $l(C(\alpha)) = m$ if and only if $\alpha_{R} = \alpha$.
    \item[(d)] A decomposition $\alpha = \alpha_{1}\alpha_{2}$ is the reduction-factorization of $\alpha$ if and only if $C(\alpha_{1})(m) = 0$ and $l(C(\alpha_{2})) = m$.
    \item[(e)] The suffix $\alpha_{R}$ is its own greatest factor with respect to $\prec_{\X}$. Consequently, $\alpha_{R}$ is the greatest factor of $\alpha$ with respect to $\prec_{\X}$.
\end{itemize}
\end{lemma}

\begin{lemma}\label{lem:C(a_R)=C(a)_R}
Let $\alpha, \beta$ be non-empty words in $\X$. Then the following assertions hold:
\begin{itemize}
\item[(a)] $C(\alpha_{R}) = C(\alpha)_{R}$, which implies $\operatorname{rf}(C(\alpha)) = (C(\alpha_{L}), C(\alpha_{R}))$.
\item[(b)] The suffix $C(\alpha_{R})$ is its own greatest factor with respect to $\prec_{\lOr}$, and hence $C(\alpha_{R})$ is the greatest factor of $C(\alpha)$ with respect to $\prec_{\lOr}$.
\end{itemize}
\end{lemma}

We now recall the notions of reducible and irreducible words.

Let $I$ be an ideal of the free algebra $\mathbb{K}\langle X \rangle$. A word $\alpha \in \langle X \rangle$ is called \textit{$I$-reducible} if $\alpha$ is the leading (greatest) word appearing in some non-zero polynomial in $I$ with respect to the reduction order $\prec_{\langle X \rangle}$. Otherwise, $\alpha$ is called \textit{$I$-irreducible}. The set of $I$-irreducible letters in $X$ is denoted by
\[
X_{I} := \{x \in X \mid x \text{ is } I\text{-irreducible}\}.
\]

Let $\mathbb{K}\langle X \rangle$ be the augmented free algebra equipped with the canonical augmentation map $\epsilon \colon \mathbb{K}\langle X \rangle \to \mathbb{K}$, where $\epsilon(X) = 0$. The set of irreducible letters in $O$ satisfies the following linear independence property:

\begin{lemma}\label{lem:irrl}
Let $(\mathbb{K}\langle X \rangle, \epsilon)$ be the augmented free algebra, and let $I \subseteq \ker \epsilon$ be a proper ideal of $\mathbb{K}\langle X \rangle$. If $y \in O$ is an $I$-irreducible letter, then
\[
y \notin \{x \in X \mid x \prec_{X} y\} \cdot \mathbb{K}\langle X \rangle + I.
\]
\end{lemma}

Using the Noetherian condition, we show that the underlying set of irreducible letters in $O$ is finite:

\begin{lemma}\label{lem:Noetherian-irrls-finite}
Let $(\mathbb{K}\langle X \rangle, \epsilon)$ be the augmented free algebra, and let $I \subseteq \ker\epsilon$ be a proper ideal of $\mathbb{K}\langle X \rangle$. If the quotient algebra $\mathbb{K}\langle X \rangle / I$ is right Noetherian, then $X_{I} \cap O$ is a finite set.
\end{lemma}

\begin{proof}
Suppose, on the contrary, that $X_{I} \cap O$ is infinite. We can then choose a sequence of distinct letters $x_{1}, x_{2}, \ldots, x_{n}, \ldots$ in $X_{I} \cap O$ such that
\[
x_{1} \prec_{O} x_{2} \prec_{O} \cdots \prec_{O} x_{n} \prec_{O} \cdots,
\]
and consider the ascending chain of right ideals of $\mathbb{K}\langle X \rangle$:
\begin{equation}\label{eq:chain}
x_{1}\mathbb{K}\langle X \rangle + I \subseteq (x_{1}\mathbb{K}\langle X \rangle + x_{2}\mathbb{K}\langle X \rangle) + I \subseteq \cdots \subseteq \sum_{i=1}^{n} x_{i}\mathbb{K}\langle X \rangle + I \subseteq \cdots.
\end{equation}
Since $X, I \subseteq \ker \epsilon$, we have $1 \notin \sum_{i=1}^{n} x_{i}\mathbb{K}\langle X \rangle + I$, which implies that each right ideal in \eqref{eq:chain} is proper. By Lemma~\ref{lem:irrl}, it follows that
\[
x_{n+1} \notin \sum_{i=1}^{n} x_{i}\mathbb{K}\langle X \rangle + I, \qquad \text{for all } n \ge 1.
\]
Hence, the chain \eqref{eq:chain} is strictly ascending. Passing to the quotient algebra $\mathbb{K}\langle X \rangle / I$, this yields an infinite strictly ascending chain of right ideals in $\mathbb{K}\langle X \rangle / I$, contradicting the hypothesis that $\mathbb{K}\langle X \rangle / I$ is right Noetherian. Therefore, $X_{I} \cap O$ must be finite.
\end{proof}

\section{Comultiplication on a Free Algebra in ${}_{H}^{H}\mathcal{YD}$}\label{sec:3}

In this section, using reduction orders and prime words, we determine the comultiplication on arbitrary non-empty words. 

For a set $X$ equipped with a degree map $t: X \rightarrow \mathbb{N}$, we write
\begin{align*}
X_{m} := \{x \in X \mid t(x) = m\}, \quad m \ge 0, \qquad X_{+} := \{x \in X \mid t(x) \ge 1\}.
\end{align*}  
The free algebra $\mathbb{K}\X$ admits a grading $\mathbb{K}\X = \bigoplus_{n=0}^{\infty} \mathbb{K}\X_{n}$, defined by $\deg(x) = t(x)$ for every letter $x \in X$. 

Throughout this section, assume that $(O, \prec_{O})$ and $(X = \bigcup_{z \in O} C_z, \prec_{X})$ are the well-ordered sets defined in Section~\ref{sec:2}, equipped with a map $t: X \rightarrow \mathbb{N}$ satisfying:
\begin{align*}
t(x) = t(z), \quad \text{for all } x \in C_z. 
\end{align*}

Given a Hopf algebra $H$ over $\mathbb{K}$, assume further that $\mathbb{K} X$ is an object of ${_{H}^{H}\mathcal{YD}}$. Then the free algebra $\mathbb{K}\X$ is an algebra in ${_{H}^{H}\mathcal{YD}}$, whose Yetter--Drinfel'd module structure is induced by $\mathbb{K} X$.

We now introduce the following comultiplication on $\mathbb{K}\X$:

\begin{definition}
With notation as above, a \textit{comultiplication} on the free algebra $\mathbb{K}\X$ in ${_{H}^{H}\mathcal{YD}}$ is an algebra map $\Delta: \mathbb{K}\X \rightarrow \mathbb{K}\X \otimes^{c} \mathbb{K}\X$.
A comultiplication $\Delta$ on $\mathbb{K}\X$ is called \textit{YD-triangular} if the following conditions hold:
\begin{itemize}
\item For each $x \in C_z$ with $z \in O$, $x \in (H^* \otimes \id)(\rho(z))$;  
\item For each $x \in C_z$ with $z \in O$, there exist $x_{(-1)} \in H$ and $x_{(0)} \in \bigcup_{d \preceq_{O} z} C_d$ such that
\begin{align*}
\rho(x) = \sum x_{(-1)} \otimes x_{(0)};
\end{align*}
\item For each $x \in X$, there exist polynomials $x' \in \mathbb{K}\X^{+}$ and letters $x'' \in X$ such that   
\begin{align*}
\Delta(x) = x \otimes 1 + 1 \otimes x + \sum x' \otimes x'',
\end{align*}
where $t(x'') < t(x)$ and $C(x'') \prec_{O} C(x)$.
\end{itemize}
\end{definition}

For each word $\omega = x_{1}\cdots x_{n} \in \X$ with letters $x_{1}, \ldots, x_{n} \in X$, write
\begin{align*}
\omega_{(0)} &:= \prod_{i=1}^{n} (x_{i})_{(0)}.
\end{align*}
Since $C(\omega_{(0)}) \preceq_{\lOr} C(\omega)$ and $|\omega_{(0)}| = |\omega|$, Lemma~\ref{lem:C(a_R)=C(a)_R} implies that
\begin{align}\label{formula:C-w_0-R}
C((\omega_{(0)})_{R}) = C(\omega_{(0)})_{R} \preceq_{\lOr} C(\omega)_{R} = C(\omega_{R}).
\end{align}

\begin{lemma}\label{lem:deg}
Let $\Delta$ be a YD-triangular comultiplication on $\mathbb{K}\X$, and let $f \in \X^{+}$. Suppose there exist polynomials $f' \in \mathbb{K}\X^{+}$ and words $f'' \in \X^{+}$ such that
\begin{align*}
\Delta(f) = f \otimes 1 + 1 \otimes f + \sum f' \otimes f''.
\end{align*}
Then $\deg(f'') < \deg(f)$ and $|f''| \le |f|$.
\end{lemma}

\begin{proof}
The comultiplication on $X$ shows that $|1| < |x''| = |x|$ for each $x \in X_{+}$, yielding $|f''| \le |f|$.
Since $0 = \deg(1) \le \deg(x'') < \deg(x)$ for each $x \in X_{+}$, it follows that $\deg(f'') < \deg(f)$.
\end{proof}

We next determine the comultiplication of prime words.

\begin{lemma}\label{lem:coproduct-omega-1}
Let $\Delta$ be a YD-triangular comultiplication on $\mathbb{K}\langle X \rangle$, and let $\omega \in \langle X \rangle_{p}$. Then there exist polynomials $\omega' \in \mathbb{K} \langle X \rangle^{+}$ and words $\omega'' \in \langle X \rangle^{+}$ such that
\begin{align*}
\Delta(\omega) = \omega \otimes 1 + 1 \otimes \omega + \sum \omega' \otimes \omega'',
\end{align*}
where $|\omega''| \le |\omega|$ and $C(\omega'') \preceq_{\langle O \rangle} C((\omega'')_{R}) \prec_{\langle O \rangle} C(\omega)$, which implies $\omega'' \preceq_{\langle X \rangle} (\omega'')_{R} \prec_{\langle X \rangle} \omega$.
\end{lemma}

\begin{proof}
We proceed by induction on the length of words. For $|\omega| = 1$, the assertion follows directly from the comultiplication formula for letters of $X$. Assume $|\omega| \ge 2$ and write $\omega = \alpha x$ with $x \in X$. Since $l(C(\alpha)) \preceq_{X} C(m_{\alpha}) \preceq_{X} C(m_{\omega}) = l(C(\omega)) = l(C(\alpha))$, it follows that $l(C(\alpha)) = C(m_{\alpha})$. Lemma~\ref{lem:properties-alpha_R}(c) then implies that $\alpha$ is a prime word.

By the inductive hypothesis applied to $\alpha$ (since $|\alpha| < |\omega|$) and the comultiplication on $X$, there exist polynomials $x', \alpha' \in \mathbb{K}\langle X \rangle^{+}$, letters $x'' \in X$, and words $\alpha'' \in \langle X \rangle^{+}$ such that
\begin{align*}
\Delta(\omega) &= \Delta(\alpha) \Delta(x) \\
&= \left(\alpha \otimes 1 + 1 \otimes \alpha + \sum \alpha' \otimes \alpha''\right)\left(x \otimes 1 + 1 \otimes x + \sum x' \otimes x''\right) \\
&= \alpha x \otimes 1 + \alpha \otimes x + \sum \alpha x' \otimes x'' \\
&\quad + \sum \alpha_{(-1)} \cdot x \otimes \alpha_{(0)} + 1 \otimes \alpha x + \sum \sum \alpha_{(-1)} \cdot x' \otimes \alpha_{(0)} x'' \\
&\quad + \sum \sum \alpha' ((\alpha'')_{(-1)} \cdot x) \otimes (\alpha'')_{(0)} + \sum \alpha' \otimes \alpha'' x \\
&\quad + \sum \sum \alpha' ((\alpha'')_{(-1)} \cdot x') \otimes (\alpha'')_{(0)} x'',
\end{align*}
where $C(x'), C(x'') \prec_{O} C(x)$, $C(\alpha'') \preceq_{\langle O \rangle} C((\alpha'')_{R}) \prec_{\langle O \rangle} C(\alpha)$, and $|\alpha''| \le |\alpha|$.

By Lemma~\ref{lem:C(a_R)=C(a)_R}, it suffices to verify that each right tensor factor $\omega''$ arising in the above sum satisfies
\[
C((\omega'')_{R}) \prec_{\langle O \rangle} C(\alpha) \quad \text{and} \quad |\omega''| \le |\omega|.
\]
Since $\omega$ is prime, Lemma~\ref{lem:C(a_R)=C(a)_R} yields $C(\alpha) \prec_{\langle O \rangle} C(\omega)$ and $C(x) \prec_{\langle O \rangle} C(\omega)$. We now check these conditions case by case using equation~\eqref{formula:C-w_0-R}, Lemma~\ref{lem:ac<bc}, and Lemma~\ref{lem:C(a_R)=C(a)_R}:

\begin{itemize}
  \item \textit{Case 1:} $\omega'' = x$. \\
  Then $C(x_{R}) = C(x) \prec_{\langle O \rangle} C(\omega)$.

  \item \textit{Case 2:} $\omega'' = x''$. \\
  Then $C((x'')_{R}) = C(x'') \prec_{\langle O \rangle} C(x) \prec_{\langle O \rangle} C(\omega)$.

  \item \textit{Case 3:} $\omega'' = \alpha_{(0)}$. \\
  Then $C((\alpha_{(0)})_{R}) \preceq_{\langle O \rangle} C(\alpha_{R}) = C(\alpha) \prec_{\langle O \rangle} C(\omega)$.

  \item \textit{Case 4:} $\omega'' = \alpha_{(0)} x''$. \\
  Then $C((\alpha_{(0)} x'')_{R}) \preceq_{\langle O \rangle} (C(\alpha) C(x''))_{R}$.
  \begin{itemize}
    \item If $m_{C(\alpha)} \succeq_{O} C(x'')$, then $(C(\alpha)C(x''))_{R} = C(\alpha)_{R} C(x'') = C(\alpha) C(x'') \prec_{\langle O \rangle} C(\alpha) C(x) = C(\omega)$.
    \item If $m_{C(\alpha)} \prec_{O} C(x'')$, then $(C(\alpha)C(x''))_{R} = C(x'') \prec_{\langle O \rangle} C(\omega)$.
  \end{itemize}

  \item \textit{Case 5:} $\omega'' = (\alpha'')_{(0)}$. \\
  Then $C(((\alpha'')_{(0)})_{R}) \preceq_{\langle O \rangle} C((\alpha'')_{R}) \prec_{\langle O \rangle} C(\alpha) \prec_{\langle O \rangle} C(\omega)$.

  \item \textit{Case 6:} $\omega'' = \alpha'' x$. \\
  Then $C((\alpha'' x)_{R}) = (C(\alpha'') C(x))_{R}$.
  \begin{itemize}
    \item If $m_{C(\alpha'')} \succeq_{O} C(x)$, then $(C(\alpha'') C(x))_{R} = C(\alpha'')_{R} C(x) \prec_{\langle O \rangle} C(\alpha) C(x) = C(\alpha x)$.
    \item If $m_{C(\alpha'')} \prec_{O} C(x)$, then $(C(\alpha'') C(x))_{R} = C(x) \prec_{\langle O \rangle} C(\alpha x)$.
  \end{itemize}

  \item \textit{Case 7:} $\omega'' = (\alpha'')_{(0)} x''$. \\
  Then $C(((\alpha'')_{(0)} x'')_{R}) \preceq_{\langle O \rangle} (C(\alpha'') C(x''))_{R}$.
  \begin{itemize}
    \item If $m_{C(\alpha'')} \succeq_{O} C(x'')$, then $(C(\alpha'') C(x''))_{R} = C(\alpha'')_{R} C(x'') = C((\alpha'')_{R}) C(x'') \prec_{\langle O \rangle} C(\alpha) C(x'') \prec_{\langle O \rangle} C(\alpha) C(x) = C(\alpha x)$.
    \item If $m_{C(\alpha'')} \prec_{O} C(x'')$, then $(C(\alpha'') C(x''))_{R} = C(x'') \prec_{\langle O \rangle} C(\alpha x)$.
  \end{itemize}

  \item \textit{Length Bounds:} \\
  Direct inspection yields $|x|, |x''|, |\alpha_{(0)}|, |\alpha_{(0)} x''|, |(\alpha'')_{(0)}|, |\alpha'' x|, |(\alpha'')_{(0)} x''| \le |\alpha x| = |\omega|$.
\end{itemize}
Therefore, every tensor factor $\omega''$ satisfies the claimed properties, completing the proof.
\end{proof}

By Lemma~\ref{lem:coproduct-omega-1}, we now determine the comultiplication structure for non-prime words.

\begin{lemma}\label{lem:coproduct-omega-2}
Suppose that $\Delta$ is a YD-triangular comultiplication on $\mathbb{K}\langle X \rangle$. Let $\omega \in \langle X \rangle \setminus \langle X \rangle_{p}$ with reduction-factorization $\mathrm{rf}(\omega) = (\omega_{L}, \omega_{R})$. Then the following assertions hold:
\begin{itemize}
  \item[(a)] There exist polynomials $(\omega_{L})', (\omega_{R})' \in \mathbb{K}\langle X \rangle^{+}$ and words $(\omega_{L})'', (\omega_{R})'' \in \langle X \rangle^{+}$ such that
  \begin{align*}
  \Delta(\omega_{L}) &= \omega_{L} \otimes 1 + 1 \otimes \omega_{L} + \sum (\omega_{L})' \otimes (\omega_{L})'', \\
  \Delta(\omega_{R}) &= \omega_{R} \otimes 1 + 1 \otimes \omega_{R} + \sum (\omega_{R})' \otimes (\omega_{R})'',
  \end{align*}
  where $C(m_{(\omega_{L})''}) \preceq_{O} C(m_{\omega_{L}}) \prec_{O} C(m_{\omega})$, $C(m_{(\omega_{R})''}) \preceq_{O} C(m_{\omega_{R}}) = C(m_{\omega})$, and $C((\omega_{R})'') \preceq_{\langle O \rangle} C(((\omega_{R})'')_{R}) \prec_{\langle O \rangle} C(\omega_{R})$.

  \item[(b)] There exist polynomials $\omega' \in \mathbb{K}\langle X \rangle^{+}$ and words $\omega'' \in \langle X \rangle^{+}$ such that
  \begin{align*}
  \Delta(\omega) = \omega \otimes 1 + \omega_{L} \otimes \omega_{R} + 1 \otimes \omega + \sum \omega' \otimes \omega'',
  \end{align*}
  where either $C((\omega'')_{R}) \prec_{\langle O \rangle} C(\omega_{R})$ holds, or $\omega'' = f\omega_{R}$ for some $f \in \langle X \rangle^{+}$ satisfying $C(m_{f}) \prec_{O} C(m_{\omega_{R}})$.
\end{itemize}
\end{lemma}

\section{Words and Irreducible Letters}\label{sec:4}

Throughout this section, we retain the notation established in Section~\ref{sec:3} and further assume that $\mathbb{K}\langle X \rangle$ is the augmented free algebra in ${}_{H}^{H}\mathcal{YD}$, equipped with the canonical augmentation $\epsilon \colon \mathbb{K}\langle X \rangle \to \mathbb{K}$ (so that $X \subseteq \ker \epsilon$). Using the comultiplication structure, we prove that every non-empty word can be expressed as a linear combination of products of irreducible letters.

We begin with the following reduction lemma for prime words, which serves as a crucial technical step toward Lemma~\ref{lem:rl-sum-smaller-irrls} and Theorem~\ref{thm:rl-sum-smaller-irrls}.

\begin{lemma}\label{lem:wr-sum-smaller-words}
Let $(\mathbb{K}\langle X \rangle, \epsilon)$ be the augmented free algebra in ${}_{H}^{H}\mathcal{YD}$, and let $I \subseteq \ker \epsilon$ be a proper ideal of $\mathbb{K}\langle X \rangle$. Suppose $\Delta$ is a YD-triangular comultiplication on $\mathbb{K}\langle X \rangle$ such that $\Delta(I) \subseteq I \otimes \mathbb{K}\langle X \rangle + \mathbb{K}\langle X \rangle \otimes I$.
Let $\beta$ be a prime word in $\langle X \rangle^{+}$, and suppose there exist words $f, g \in \langle X \rangle^{+}$ and scalars $k_{f}, k_{g} \in \mathbb{K}$ such that
\begin{equation}\label{formula:wr-fwr-g-I}
p_{\beta} := \beta + \sum_{f} k_{f} f \beta + \sum_{g} k_{g} g \in I, \qquad \text{with } C(m_{f}) \prec_{O} C(m_{\beta}) \text{ and } g_{R} \prec_{\langle X \rangle} \beta.
\end{equation}
Then there exist words $h \in \langle X \rangle^{+}$ and scalars $k_{h} \in \mathbb{K}$ such that
\begin{equation}\label{formula:wr-h-I}
\beta \in \sum_{h} k_{h} h + I, \qquad \text{with } h_{R} \prec_{\langle X \rangle} \beta.
\end{equation}
\end{lemma}

\begin{proof}
Define 
\[
W_{1} := \{f \mid p_{\beta} \in I \text{ and } k_{f} \neq 0\} \quad \text{and} \quad d_{1} := \max\{\deg(f) \mid f \in W_{1}\}.
\]
If $W_{1} = \emptyset$, the assertion holds trivially. Assume that $W_{1} \neq \emptyset$. By Lemma~\ref{lem:coproduct-omega-2}, there exist polynomials $f', g', \beta' \in \mathbb{K}\langle X \rangle^{+}$ and words $f'', g'', \beta'' \in \langle X \rangle^{+}$ such that
\begin{align*}
&\Delta\left(\beta + \sum_{f \in W_{1}} k_{f} f \beta + \sum_{g} k_{g} g\right) \\
&= \beta \otimes 1 + \sum \beta' \otimes \beta'' + 1 \otimes \beta \\
&\quad + \sum_{f \in W_{1}} k_{f} \left(f \otimes 1 + 1 \otimes f + \sum f' \otimes f''\right) \left(\beta \otimes 1 + 1 \otimes \beta + \sum \beta' \otimes \beta''\right) \\
&\quad + \sum_{g} k_{g} \left(g \otimes 1 + 1 \otimes g + \sum g' \otimes g''\right) \\
&= \beta \otimes 1 + 1 \otimes \beta + \sum \beta' \otimes \beta'' \\
&\quad + \sum_{f \in W_{1}} k_{f} \Big(f \beta \otimes 1 + f \otimes \beta + \sum f \beta' \otimes \beta'' \\
&\qquad + \sum f_{(-1)} \cdot \beta \otimes f_{(0)} + 1 \otimes f \beta + \sum f_{(-1)} \cdot \beta' \otimes f_{(0)} \beta'' \\
&\qquad + \sum f' ((f'')_{(-1)} \cdot \beta) \otimes (f'')_{(0)} + \sum f' \otimes f'' \beta + \sum f' ((f'')_{(-1)} \cdot \beta') \otimes (f'')_{(0)} \beta''\Big) \\
&\quad + \sum_{g} k_{g} \left(g \otimes 1 + 1 \otimes g + \sum g' \otimes g''\right),
\end{align*}
where $C(m_{f''}) \preceq_{O} C(m_{f}) \prec_{O} C(m_{\beta})$, $C((\beta'')_{R}) \prec_{\langle O \rangle} C(\beta)$, $(g'')_{R} \preceq_{\langle X \rangle} g_{R} \prec_{\langle X \rangle} \beta$, and $\deg(f'') < d_{1}$ by Lemma~\ref{lem:deg}.

Set $p_{1} := \sum_{f \in W_{1}} k_{f} f$. If $p_{1} \in I$, the result follows directly. Assume therefore that $p_{1} \notin I$. Evaluating under the augmentation $\epsilon$, the sum $\mathbb{K}p_{1} \oplus \mathbb{K}1 \oplus I$ is direct. Decomposing the free algebra as $\mathbb{K}\langle X \rangle = \mathbb{K}p_{1} \oplus \mathbb{K}1 \oplus I \oplus U$ for some subspace $U \subseteq \mathbb{K}\langle X \rangle$, we define a linear functional $\phi \colon \mathbb{K}\langle X \rangle \to \mathbb{K}$ by setting
\[
\phi(p_{1}) = 1, \quad \phi(1) = 0, \quad \text{and} \quad \phi(I) = \phi(U) = 0.
\]
We observe that:
\begin{itemize}
  \item $(\beta'')_{R} \prec_{\langle X \rangle} \beta$, because $C((\beta'')_{R}) \prec_{\langle O \rangle} C(\beta)$;
  \item $(f_{(0)})_{R}, (f_{(0)} \beta'')_{R}, ((f'')_{(0)})_{R}, ((f'')_{(0)} \beta'')_{R} \prec_{\langle X \rangle} \beta$, because $C(m_{(f'')_{(0)}}), C(m_{f''}), C(m_{f_{(0)}}) \prec_{O} C(m_{f}) \prec_{O} C(m_{\beta})$, and $(\beta'')_{R} \prec_{\langle X \rangle} \beta$;
  \item $(g'')_{R} \preceq_{\langle X \rangle} g_{R} \prec_{\langle X \rangle} \beta$.
\end{itemize}

Applying the linear map $\phi \otimes \operatorname{id}$ to the containment relation
\[
\Delta\left(\beta + \sum_{f \in W_{1}} k_{f} f \beta + \sum_{g} k_{g} g\right) \in I \otimes \mathbb{K}\langle X \rangle + \mathbb{K}\langle X \rangle \otimes I,
\]
we obtain words $g_{1} \in \langle X \rangle^{+}$ and scalars $k_{g_{1}} \in \mathbb{K}$ such that
\begin{equation}\label{formula:wr-f''wr-g1-I}
\beta + \sum_{f \in W_{1}} k_{f} \phi(f') f'' \beta + \sum k_{g_{1}} g_{1} \in I,
\end{equation}
where $C(m_{f''}) \prec_{O} C(m_{\beta})$, $(g_{1})_{R} \prec_{\langle X \rangle} \beta$, and $\deg(f'') < d_{1}$.

Set $d_{2} := \max\{\deg(f'') \mid f \in W_{1}\}$. We now prove \eqref{formula:wr-h-I} by induction on $d_{1}$. If $d_{1} = 1$, then $\deg(f) \le 1$ for all $f \in W_{1}$, which implies $f \in X$. Since $\Delta(f) = f \otimes 1 + 1 \otimes f$, equation~\eqref{formula:wr-f''wr-g1-I} simplifies directly to
\[
\beta + \sum k_{g_{1}} g_{1} \in I, \qquad \text{with } (g_{1})_{R} \prec_{\langle X \rangle} \beta.
\]
For $d_{1} > 1$, since $d_{2} < d_{1}$, the inductive hypothesis applied to equation~\eqref{formula:wr-f''wr-g1-I} yields the required form~\eqref{formula:wr-h-I}, completing the proof.
\end{proof}

We next show that any reducible letter can be written as a linear combination of products of irreducible letters of strictly smaller order.

\begin{lemma}\label{lem:rl-sum-smaller-irrls}
Let $(\mathbb{K}\langle X \rangle, \epsilon)$ be the augmented free algebra in ${}_{H}^{H}\mathcal{YD}$, and let $I \subseteq \ker \epsilon$ be a proper ideal of $\mathbb{K}\langle X \rangle$.
Suppose $\Delta$ is a YD-triangular comultiplication on $\mathbb{K}\langle X \rangle$ such that $\Delta(I) \subseteq I \otimes \mathbb{K}\langle X \rangle + \mathbb{K}\langle X \rangle \otimes I$.
If $x$ is an $I$-reducible letter in $X$, then there exist letters $x_{1}, \ldots, x_{n} \in X$ and scalars $k_{x_{1}\cdots x_{n}} \in \mathbb{K}$ such that
\begin{equation}
x \in \sum k_{x_{1}\cdots x_{n}} x_{1} \cdots x_{n} + I, \qquad \text{with } x_{1}, \ldots, x_{n} \prec_{X} x.
\end{equation}
\end{lemma}

\begin{proof}
If $x \in I$, the assertion holds trivially. Assume therefore that $x \notin I$. We proceed by case analysis on $x$:

\paragraph{\textit{Case 1}: $x \in O$.}
Evaluating under the augmentation $\epsilon$, there exist words $\omega \in \langle X \rangle^{+}$ and scalars $k_{\omega} \in \mathbb{K}$ such that
\[
x \in \sum k_{\omega} \omega + I, \qquad \text{with } \omega \prec_{\langle X \rangle} x,
\]
or equivalently, by Corollary~\ref{cor:well-orders-compatibility-2},
\[
x \in \sum k_{\omega} \omega + I, \qquad \text{with } C(l(\omega)) \prec_{O} x.
\]
Define the support set $W$ and its associated set of prime factors $W_{R}$ by
\[
W := \{\omega \mid x \in \sum k_{\omega} \omega + I \text{ and } k_{\omega} \neq 0\}, \qquad W_{R} := \{\omega_{R} \mid \omega \in W\}.
\]
Let $\beta$ denote the greatest word of $W_{R}$ with respect to $\prec_{\langle X \rangle}$. We can then decompose $x$ modulo $I$ as
\[
x \in \sum_{\alpha\beta \in W} k_{\alpha\beta} \alpha\beta + \sum_{\substack{\omega \in W \\ \omega_{R} \prec_{\langle X \rangle} \beta}} k_{\omega} \omega + I, \qquad \text{with } C(l(\alpha\beta)), C(l(\omega)) \prec_{O} x.
\]
Without loss of generality, we may assume that the leading sum $\sum_{\alpha\beta \in W} k_{\alpha\beta} \alpha\beta \notin I$.

Since $(\langle X \rangle, \prec_{\langle X \rangle})$ is well-ordered, we proceed by induction on $\beta$ with respect to $\prec_{\langle X \rangle}$.

\noindent
\textit{Subcase 1.1: $\beta \prec_{\langle X \rangle} x$.} For each word $\omega \in W$ and each letter $b$ appearing in $\omega$, Lemma~\ref{lem:properties-alpha_R}(e) yields $b \preceq_{\langle X \rangle} \omega_{R} \preceq_{\langle X \rangle} \beta \prec_{\langle X \rangle} x$. By Lemma~\ref{lem:compatibility-orders}, this implies $b \prec_{X} x$, completing this subcase.

\noindent
\textit{Subcase 1.2: $\beta \succeq_{\langle X \rangle} x$.} By Corollary~\ref{cor:well-orders-compatibility-2}, this is equivalent to $C(l(\beta)) \succeq_{O} x$. This implies $\alpha \neq 1$; otherwise, $C(l(\beta)) \prec_{O} x$, a contradiction.

Let $a$ denote the least non-empty word in $\langle X \rangle$, which is also the least letter in $X$ with respect to $\prec_{X}$ by Lemma~\ref{lem:x-least}(c). If $\beta = a$, then for each $\omega \in W$, we have $\omega \in \langle a \rangle$ and $a = C(l(\omega)) \prec_{O} x$, which settles this situation.

Assume now that $\beta \succ_{\langle X \rangle} a$. By Lemma~\ref{lem:coproduct-omega-2} and the comultiplication formula on $X$, there exist polynomials $x', (\alpha\beta)', \omega' \in \mathbb{K}\langle X \rangle^{+}$, letters $x'' \in X$, and words $(\alpha\beta)'', \omega'' \in \langle X \rangle^{+}$ such that
\begin{align*}
&\Delta\left(\sum_{\alpha\beta \in W} k_{\alpha\beta} \alpha\beta + \sum_{\substack{\omega \in W \\ \omega_{R} \prec_{\langle X \rangle} \beta}} k_{\omega} \omega - x\right) \\
&= \sum_{\alpha\beta \in W} k_{\alpha\beta} \left(\alpha\beta \otimes 1 + \alpha \otimes \beta + 1 \otimes \alpha\beta + \sum (\alpha\beta)' \otimes (\alpha\beta)''\right) \\
&\quad + \sum_{\omega_{R} \prec_{\langle X \rangle} \beta} k_{\omega} \left(\omega \otimes 1 + 1 \otimes \omega + \sum \omega' \otimes \omega''\right) - \left(x \otimes 1 + 1 \otimes x + \sum x' \otimes x''\right),
\end{align*}
where $x'' \prec_{X} x$, and for each $(\alpha\beta)''$, either $((\alpha\beta)'')_{R} \prec_{\langle X \rangle} \beta$ or $(\alpha\beta)'' = u\beta$ for some $u \in \langle X \rangle^{+}$ satisfying $C(m_{u}) \prec_{O} C(m_{\beta})$. Furthermore, $(\omega'')_{R} \preceq_{\langle X \rangle} \omega_{R} \prec_{\langle X \rangle} \beta$.

Set $p := \sum_{\alpha\beta \in W} k_{\alpha\beta} \alpha$. Since $\sum_{\alpha\beta \in W} k_{\alpha\beta} \alpha\beta \notin I$, the sum $\mathbb{K}p \oplus I$ is direct. Evaluating under $\epsilon$, the sum $\mathbb{K}p \oplus I \oplus \mathbb{K}1$ is also direct. Decomposing the free algebra as $\mathbb{K}\langle X \rangle = \mathbb{K}p \oplus I \oplus \mathbb{K}1 \oplus V$ for some subspace $V \subseteq \mathbb{K}\langle X \rangle$, we define a linear functional $\varphi \colon \mathbb{K}\langle X \rangle \to \mathbb{K}$ by
\[
\varphi(p) = 1, \quad \varphi(1) = 0, \quad \text{and} \quad \varphi(I) = \varphi(V) = 0.
\]
Applying the map $\varphi \otimes \operatorname{id}$ to the containment relation
\[
\Delta\left(\sum_{\alpha\beta \in W} k_{\alpha\beta} \alpha\beta + \sum_{\substack{\omega \in W \\ \omega_{R} \prec_{\langle X \rangle} \beta}} k_{\omega} \omega - x\right) \in I \otimes \mathbb{K}\langle X \rangle + \mathbb{K}\langle X \rangle \otimes I,
\]
we obtain words $f, g \in \langle X \rangle^{+}$ and scalars $k_{f}, k_{g} \in \mathbb{K}$ such that
\[
\beta + \sum_{f} k_{f} f\beta + \sum_{g} k_{g} g \in I, \qquad \text{with } C(m_{f}) \prec_{O} C(\beta) \text{ and } g_{R} \prec_{\langle X \rangle} \beta.
\]
By Lemma~\ref{lem:wr-sum-smaller-words}, there exist words $h \in \langle X \rangle^{+}$ and scalars $k_{h} \in \mathbb{K}$ such that
\[
\beta \in \sum_{h} k_{h} h + I, \qquad \text{with } h_{R} \prec_{\langle X \rangle} \beta.
\]
Substituting this back into the expansion for $x$ yields
\[
x \in \sum_{\alpha\beta \in W} \sum_{h} k_{\alpha\beta} k_{h} \alpha h + \sum_{\substack{\omega \in W \\ \omega_{R} \prec_{\langle X \rangle} \beta}} k_{\omega} \omega + I.
\]
Note that $(\alpha h)_{R} \prec_{\langle X \rangle} \beta$ because $h_{R} \prec_{\langle X \rangle} \beta$ and $C(m_{\alpha}) \prec_{O} C(m_{\beta})$. By the inductive hypothesis (since $(\alpha h)_{R}, \omega_{R} \prec_{\langle X \rangle} \beta$), there exist letters $x_{1}, \ldots, x_{n} \in X$ and scalars $k_{x_{1}\cdots x_{n}} \in \mathbb{K}$ such that
\[
x \in \sum k_{x_{1}\cdots x_{n}} x_{1} \cdots x_{n} + I, \qquad \text{with } x_{1}, \ldots, x_{n} \prec_{X} x.
\]

\paragraph{\textit{Case 2}: $x \in C_{z} \setminus O$ with $z \in O$.}
By evaluation under $\epsilon$, we have $x \in \sum k_{\omega} \omega + I$ for some words $\omega \in \langle X \rangle^{+}$ satisfying $\omega \prec_{\langle X \rangle} x$ and scalars $k_{\omega} \in \mathbb{K}$. By the definition of $\prec_{\langle X \rangle}$, there exist letters $z_{1}, \ldots, z_{s} \in C_{x} \setminus O$ such that
\[
x \in k_{z} z + \sum_{1 \le i \le s} k_{z_{i}} z_{i} + \sum k_{\omega} \omega + I, \qquad \text{with } \omega \prec_{\langle X \rangle} z.
\]
An analogous reduction argument to Case 1 applies here by tracking the comultiplication $\Delta(x - k_{z} z - \sum_{1 \le i \le s} k_{z_{i}} z_{i})$ instead of $\Delta(x)$. It follows that there exist letters $y_{1}, \ldots, y_{m} \in X$ and scalars $k_{y_{1}\cdots y_{m}} \in \mathbb{K}$ such that
\[
x \in k_{z} z + \sum_{1 \le i \le s} k_{z_{i}} z_{i} + \sum k_{y_{1} \cdots y_{m}} y_{1} \cdots y_{m} + I,
\]
where $y_{1}, \ldots, y_{m} \prec_{X} z \prec_{X} z_{1} \prec_{X} \cdots \prec_{X} z_{s} \prec_{X} x$. This completes the proof.
\end{proof}

We show that the original element of every irreducible letter is also irreducible.

\begin{lemma}\label{lem:rl-sum-smaller-irrls-2}
Let $(\mathbb{K}\langle X \rangle, \epsilon)$ be the augmented free algebra in ${}_{H}^{H}\mathcal{YD}$, and let $I \subseteq \ker\epsilon$ be a proper ideal of $\mathbb{K}\langle X \rangle$. Suppose $\Delta$ is a YD-triangular comultiplication on $\mathbb{K}\langle X \rangle$ such that $\Delta(I) \subseteq I \otimes \mathbb{K}\langle X \rangle + \mathbb{K}\langle X \rangle \otimes I$.
If $x \in O$ is $I$-reducible, then every letter $y \in C_{x}$ is $I$-reducible.
\end{lemma}

\begin{proof}
By Lemma~\ref{lem:rl-sum-smaller-irrls}, there exist letters $x_{1}, \ldots, x_{n} \in X$ and scalars $k_{x_{1}\cdots x_{n}} \in \mathbb{K}$ such that
\[
x \in \sum k_{x_{1}\cdots x_{n}} x_{1} \cdots x_{n} + I, \qquad \text{with } x_{1}, \ldots, x_{n} \prec_{X} x.
\]
By the structural properties of $X$, there exists a linear functional $h^{*} \in H^{*}$ such that $y = (h^{*} \otimes \operatorname{id})(\rho(x))$. Since the coaction $\rho$ is an algebra map, applying $h^{*} \otimes \operatorname{id}$ to the reduction of $x$ yields letters $y_{1}, \ldots, y_{n} \in X$ and scalars $k_{y_{1}\cdots y_{n}} \in \mathbb{K}$ such that
\[
y \in \sum k_{y_{1}\cdots y_{n}} y_{1} \cdots y_{n} + I,
\]
where $y_{1}, \ldots, y_{n} \prec_{X} x \preceq_{X} y$. Consequently, $y$ can be expressed modulo $I$ in terms of strictly smaller terms, which implies that $y$ is $I$-reducible.
\end{proof}

We now prove that every non-empty word can be expressed modulo $I$ as a linear combination of products of irreducible letters.

\begin{theorem}\label{thm:rl-sum-smaller-irrls}
Let $(\mathbb{K}\langle X \rangle, \epsilon)$ be the augmented free algebra in ${}_{H}^{H}\mathcal{YD}$, and let $I \subseteq \ker\epsilon$ be a proper ideal of $\mathbb{K}\langle X \rangle$. Suppose $\Delta$ is a YD-triangular comultiplication on $\mathbb{K}\langle X \rangle$ such that $\Delta(I) \subseteq I \otimes \mathbb{K}\langle X \rangle + \mathbb{K}\langle X \rangle \otimes I$.
If $x \in X$ is an $I$-reducible letter, then there exist $I$-irreducible letters $x_{1}, \ldots, x_{n} \in X_{I}$ and scalars $k_{x_{1}\cdots x_{n}} \in \mathbb{K}$ such that
\begin{equation}
x \in \sum k_{x_{1}\cdots x_{n}} x_{1} \cdots x_{n} + I, \qquad \text{with } x_{1}, \ldots, x_{n} \prec_{X} x.
\end{equation}
Consequently, every non-empty word in $\langle X \rangle$ can be expressed modulo $I$ as a linear combination of products of $I$-irreducible letters.
\end{theorem}

\begin{proof}
We proceed by induction on $x \in X$ with respect to the well-order $\prec_{X}$.
By Lemma~\ref{lem:rl-sum-smaller-irrls}, there exist letters $x_{1}, \ldots, x_{n} \in X$ and scalars $k_{x_{1}\cdots x_{n}} \in \mathbb{K}$ such that
\[
x \in \sum k_{x_{1}\cdots x_{n}} x_{1} \cdots x_{n} + I, \qquad \text{with } x_{1}, \ldots, x_{n} \prec_{X} x.
\]
By the inductive hypothesis (since $x_{1}, \ldots, x_{n} \prec_{X} x$), each $x_{i}$ can be written modulo $I$ as a linear combination of products of $I$-irreducible letters strictly smaller than $x_i$. Substituting these expansions into the expression for $x$ yields
\[
x \in \sum k_{y_{1}\cdots y_{m}} y_{1} \cdots y_{m} + I, \qquad \text{with } y_{1}, \ldots, y_{m} \in X_{I} \text{ and } y_{1}, \ldots, y_{m} \prec_{X} x.
\]
Since $I$-irreducible letters trivially satisfy this property, it follows that every letter in $X$ can be expressed modulo $I$ as a linear combination of products of $I$-irreducible letters. Extending this multiplicatively to all words, we conclude that every non-empty word in $\langle X \rangle$ shares this property, completing the proof.
\end{proof}

\section{Noetherian Connected Prebraided Hopf Algebras are Affine}\label{sec:5}

In this section, by applying the reduction machinery developed in Section~\ref{sec:4} to connected prebraided Hopf algebras, we prove that any right Noetherian connected prebraided Hopf algebra is affine.

\begin{theorem}\label{thm:Noetherian-connected-braided-HA-affine}
Let $H$ be a Hopf algebra over $\mathbb{K}$, and let $B$ be a connected Hopf algebra in ${}_{H}^{H}\mathcal{YD}$. Then the following assertions hold:
\begin{itemize}
  \item[(a)] There exists a set $X$ constructed as in Theorem~\ref{thm:rl-sum-smaller-irrls} such that $\pi(X_{I})$ generates $B$ as an algebra, where $\pi \colon \mathbb{K}\langle X \rangle \to B$ is the canonical algebra projection and $I := \ker \pi$.
  \item[(b)] If $B$ is right Noetherian, then $B$ is affine.
\end{itemize}
\end{theorem}

\begin{proof}
By Lemma~\ref{lem:B_(n)-algebra-filtration}, the coradical filtration $\{B_{(n)}\}_{n \ge 0}$ is a bialgebra filtration of $B$ in ${}_{H}^{H}\mathcal{YD}$. Define a map $t_{B} \colon B \to \mathbb{N}$ by 
\[
t_{B}(d) := \min\{n \mid d \in B_{(n)}\}.
\]
We may choose a $\mathbb{K}$-basis $E$ of $\ker \epsilon_{B}$ satisfying the following property: for each $d \in E$, there exist elements $d', d'' \in E$ and scalars $k_{d',d''} \in \mathbb{K}$ such that
\[
\Delta_{B}(d) = d \otimes 1 + 1 \otimes d + \sum k_{d',d''} d' \otimes d'', \qquad \text{with } t_B(d'), t_{B}(d'') < t_{B}(d).
\]

We define an order $<_{E}$ on $E$ as follows. For each $n \ge 1$, fix a well-ordering $<_{n}$ on the subset $E_{n} := \{d \in E \mid t_{B}(d) = n\}$. Then for any $x, y \in E$, set
\begin{align*} 
x <_{E} y \Longleftrightarrow 
\begin{cases}
t_{B}(x) < t_{B}(y), \text{ or} \\
t_{B}(x) = t_{B}(y) \text{ and } x <_{t_{B}(x)} y.
\end{cases}
\end{align*}
One easily checks that $(E, <_{E})$ is a well-ordered set.

For each $d \in E$, we fix a tensor decomposition of rank $r = \operatorname{rank}(\rho_B(d))$ such that 
\[
\rho_B(d) = c_1 \otimes d + \sum_{i=2}^r c_i \otimes d_i,
\]
with $c_{1} \in H$, $c_{i} \in \ker \epsilon_H$, and $d_{i} \in \ker \epsilon_B$ satisfying $t_B(d_i) = t_B(d)$. We then define a finite set $D_d := \{d, d_2, \ldots, d_r\}$, equipped with a well-ordering $<_{D_d}$ such that $d$ is the minimal element in $D_d$.

We construct the basis components $C_d \subseteq D_d$ for each $d \in E$ by transfinite induction along $(E, <_E)$. For each $d \in E$, we select elements from $D_d$ that are linearly independent of all previously chosen vectors:
\[
v \in C_d \Longleftrightarrow v \notin \operatorname{span}\left( \bigcup_{a <_E d} C_a \cup \{u \in D_d \mid u <_{D_d} v\} \right).
\]
We then define an index set $O$ and a disjoint union $X$ by
\[
O := \{d \in E \mid C_d \neq \emptyset\}, \qquad X := \bigcup_{d \in O} C_d.
\]
By the linearity of $\rho_B$ and the minimality of tensor rank, for each $d \in E$, $C_d \neq \emptyset$ if and only if $d \notin \operatorname{span}\left(\bigcup_{a <_E d} C_a\right)$, in which case $d$ is the minimal element of $C_d$ with respect to $<_{D_d}$. Thus, $O = \{d \in E \mid d \in C_d\}$, and $X = \bigcup_{d \in O} C_d$ forms a $\mathbb{K}$-basis of $\ker \epsilon_{B}$. Moreover, $t_{B}(x) = t_{B}(d)$ for any $d \in O$ and $x \in C_{d}$.

For convenience, set $(O, \prec_{O}) := (O, <_{E})$, $(C_{d}, \prec_{C_{d}}) := (C_{d}, <_{D_{d}})$, and $t := t_{B}|_{X}$. It follows that $O, X \subseteq \ker\epsilon_{B}$, and $(O, \prec_{O})$ is well-ordered. The order $\prec_{X}$ on $X$ induced by $\prec_{O}$ is given for $x, y \in X$ by
\begin{align*}
x \prec_{X} y \Longleftrightarrow 
\begin{cases}
t(x) < t(y), \text{ or} \\
t(x) = t(y) \text{ and }
\begin{cases}
C(x) <_{t(x)} C(y), \text{ or} \\
C(x) = C(y) = d \text{ and } x \prec_{C_{d}} y.
\end{cases}
\end{cases}
\end{align*}
By Lemma~\ref{lem:well-orders}, $(X, \prec_{X})$, $(\langle O \rangle, \prec_{\langle O \rangle})$, and $(\langle X \rangle, \prec_{\langle X \rangle})$ are well-ordered.

Let $f_{X} \colon X \to B$ be the inclusion map. Then $f_{X}$ induces the canonical algebra projection $\pi \colon \mathbb{K}\langle X \rangle \to B$, and the free algebra $\mathbb{K}\langle X \rangle \in {}_{H}^{H}\mathcal{YD}$ generated by $\mathbb{K}X \in {}_{H}^{H}\mathcal{YD}$ admits a comultiplication $\Delta$ satisfying $\Delta_{B} \circ \pi = (\pi \otimes \pi) \circ \Delta$ and the following conditions:
\begin{itemize}
  \item for each $x \in C_z$ with $z \in O$, $x \in (H^* \otimes \id)(\rho(z))$;  
  \item for each $x \in C_{z}$ with $z \in O$, there exist $x_{(-1)} \in H$ and $x_{(0)} \in \bigcup_{d \preceq_{O} z} C_{d}$ such that $\rho(x) = \sum x_{(-1)} \otimes x_{(0)}$;
  \item $\Delta_{B}$ lifts to an algebra map $\Delta \colon \mathbb{K}\langle X \rangle \to \mathbb{K}\langle X \rangle \otimes^c \mathbb{K}\langle X \rangle$ satisfying, for each $x \in X$,
  \[
  \Delta(x) = x \otimes 1 + 1 \otimes x + \sum x' \otimes x'',
  \]
  with $x' \in \mathbb{K}\langle X \rangle^{+}$ and $x'' \in X$ satisfying $t(x'') < t(x)$ and $C(x'') \prec_{O} C(x)$.
\end{itemize}
Setting $I := \ker \pi$, $\Delta$ is a YD-triangular comultiplication such that $\Delta(I) \subseteq I \otimes \mathbb{K}\langle X \rangle + \mathbb{K}\langle X \rangle \otimes I$. Similarly, the counit $\epsilon_{B}$ lifts to an algebra map $\epsilon \colon \mathbb{K}\langle X \rangle \to \mathbb{K}$ with $X, I \subseteq \ker \epsilon$.

For Part (b), Lemma~\ref{lem:rl-sum-smaller-irrls-2} yields $X_{I} = \bigcup_{z \in O} (X_{I} \cap C_{z}) = \bigcup_{z \in X_{I} \cap O} (X_{I} \cap C_{z})$. Since $B$ is right Noetherian, Lemma~\ref{lem:Noetherian-irrls-finite} implies that the set $X_{I} \cap O$ is finite. Because each $C_{z}$ is a finite set, $X_{I}$ is finite. Applying Part (a), we conclude that $B$ is affine.
\end{proof}

\begin{remark}\label{rem:mirror-and-subalgebras}
(a) Within the framework of mirror reduction orders, a parallel argument shows that any left Noetherian connected Hopf algebra in ${}_{H}^{H}\mathcal{YD}$ is affine.

\smallskip
\noindent
(b) Suppose $H$ is a Hopf algebra with bijective antipode. Let $B$ be a right Noetherian connected Hopf algebra in ${}_{H}^{H}\mathcal{YD}$, and let $S$ be any left coideal subalgebra of $B$. By \cite[Proposition~3.6]{AAGMV2014}, $B$ is a free left and right $S$-module, and hence $S$ is right Noetherian. Since the proof of Theorem~\ref{thm:Noetherian-connected-braided-HA-affine} only involves the right tensor components of the comultiplication and of the left $H$-coaction, a similar argument shows that $S$ is also affine. In other words, any left coideal subalgebra of a right Noetherian connected Hopf algebra in ${}_{H}^{H}\mathcal{YD}$ is affine.
\end{remark}

We now show that the braided component of a Noetherian Radford biproduct is affine under a natural coradical hypothesis:

\begin{corollary}\label{cor:Noetherian-Radford-biproduct-affine}
Let $A$ be a bialgebra and $H$ a Hopf algebra over $\mathbb{K}$. Suppose that $\pi \colon A \to H$ and $\jmath \colon H \to A$ are bialgebra maps such that $\pi \circ \jmath = \operatorname{id}_{H}$. If $A$ is right Noetherian and $A_{(0)} \subseteq \jmath(H)$, then $B := A^{\mathrm{co}\,\pi}$ is connected and affine. Consequently, $A$ is affine if and only if $H$ is affine.
\end{corollary}

\begin{proof}
By Lemmas~\ref{lem:Radfordbiproduct} and \ref{lem:B-Noetherian}, and Theorem~\ref{thm:Noetherian-connected-braided-HA-affine}, it suffices to show that $B$ is connected. Since the idempotent projection $\Pi \colon A \to B$ is a surjective coalgebra map, it follows from \cite[Lemma~4.1.3(b)]{Ra2012} that $\{\Pi(A_{(n)})\}_{n \ge 0}$ forms a coalgebra filtration of $B$. Applying \cite[Proposition~4.12]{Ra2012}, we obtain
\[
B_{(0)} \subseteq \Pi(A_{(0)}) \subseteq \Pi(\jmath(H)).
\]
By \cite[Equation~(11.32)]{Ra2012}, $\Pi(\jmath(H)) = \epsilon_H(H) 1 = \mathbb{K} 1$, which forces $B_{(0)} = \mathbb{K} 1$. Thus, $B$ is connected, completing the proof.
\end{proof}

\begin{remark}
Corollary~\ref{cor:Noetherian-Radford-biproduct-affine} generalizes \cite[Theorem~2.2]{JZ2025-1}, which established that the braided component of a Noetherian graded Hopf algebra is affine.
\end{remark}

\bibliographystyle{plain}

\end{document}